\documentclass[11pt]{article}

\usepackage[english]{babel}
\usepackage[utf8]{inputenc}
\usepackage{amsmath,amsfonts,amssymb,amsthm,mathrsfs}
\usepackage{minitoc}
\usepackage[normalem]{ulem}
 \usepackage[colorlinks=true]{hyperref}
\renewcommand{\leq}{\leqslant}
\renewcommand{\geq}{\geqslant}
\usepackage[mono=false]{libertine}
\usepackage[cal=euler, calscaled=.95, scr=boondoxo]{mathalfa}
\useosf
\usepackage{mathpazo}
\usepackage{mathtools}

\usepackage{wasysym}

\usepackage[dvipsnames, svgnames]{xcolor}
\usepackage[a4paper,vmargin={2cm,2cm},hmargin={2cm,2cm}]{geometry}
\usepackage[font=sf, labelfont={sf,bf}, margin=1cm]{caption}
\usepackage{graphicx,graphics}
\usepackage{epsfig}
\usepackage{latexsym}
\usepackage{stmaryrd}
\usepackage{ae,aecompl}

\usepackage[english]{babel}

\usepackage{pstricks}
\usepackage{enumerate}
\usepackage{tikz}
\usepackage{fontawesome}
\usetikzlibrary{arrows,automata}

\DeclareFixedFont{\beaupetit}{T1}{ftp}{b}{n}{2cm}

\newtheorem{theorem}{Theorem}[]

\newtheorem{proposition}[theorem]{Proposition}

\newtheorem{lemma}[theorem]{Lemma}
\newtheorem{corollary}[theorem]{Corollary}

\theoremstyle{definition}

\newtheorem*{remark}{Remark}

\usepackage{todonotes}

\title{\textsc{ Erd{\H{o}}s--R\'enyi Poissonized}}
\author{Nicolas \textsc{Curien}\thanks{Universit\'e Paris-Saclay.\hfill  \href{mailto:nicolas.curien@gmail.com}{\texttt{nicolas.curien@gmail.com}}}
}
\date{}
\begin{document}
\maketitle 

\vspace{-0.5cm}
\begin{abstract} We introduce a variant of the Erd{\H{o}}s--R\'enyi random graph where the number of vertices is random and follows a Poisson law. A very simple Markov property of the  model entails that the Lukasiewicz exploration is made of \emph{independent} Poisson increments. Using a vanilla Poisson counting process,  this enables us to give very short proofs of classical results such as the phase transition for the giant component or the connectedness for the standard Erd{\H{o}}s--R\'enyi model. \end{abstract}

\section{The $ \mathrm{G_{Poi}}(\alpha,p)$ model and its exploration}
\begin{figure}[!h]
  \begin{center}
  \includegraphics[width=12cm]{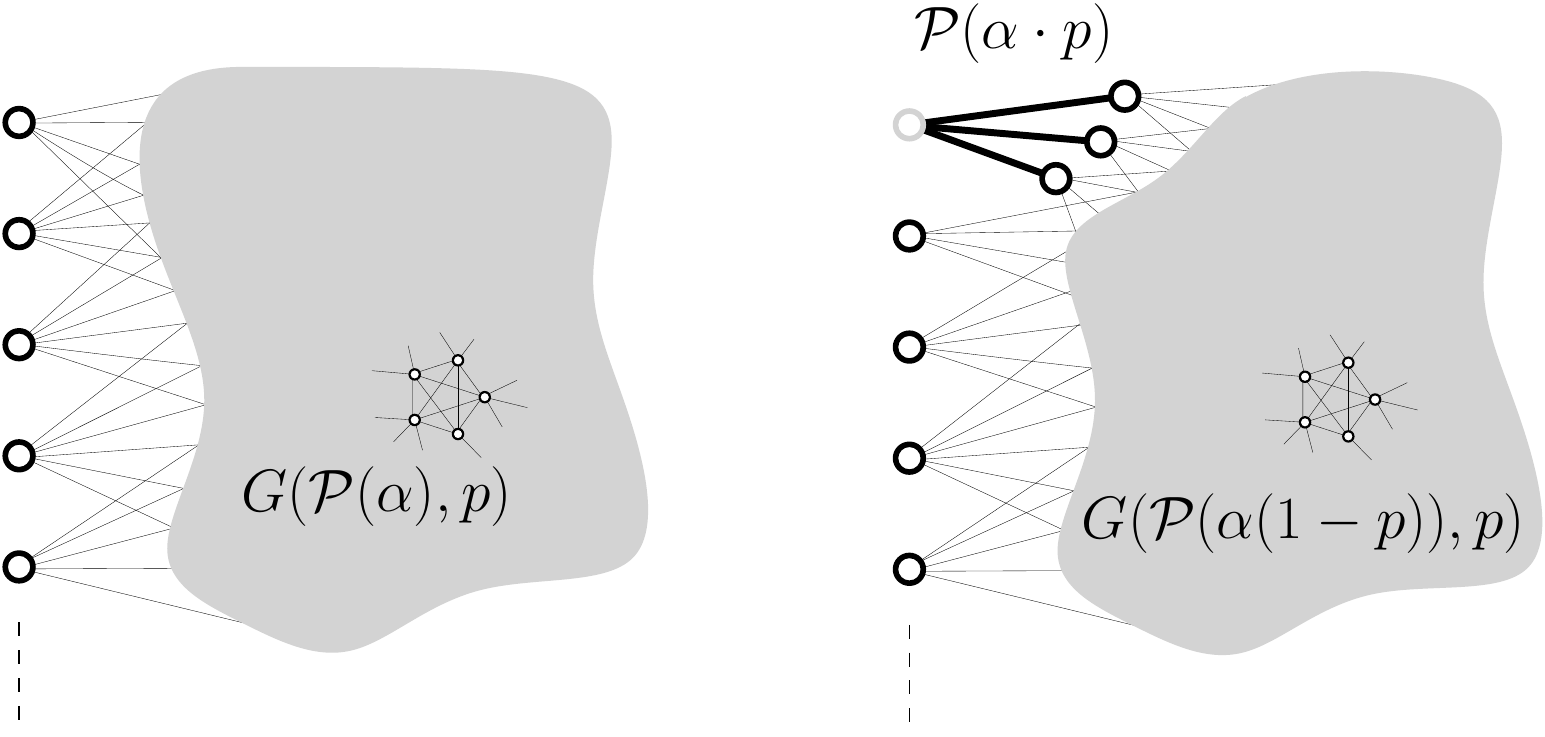}
  \caption{Discovering the neighborhood of a vertex in the stack in the Poissonized version of the Erd{\H{o}}s-R\'enyi random graph. The stack is made of the white vertices on the left part while the core is represented by the gray part. \label{fig:ERP}}
  \end{center}
  \end{figure}

Fix $\alpha >0$ and let $N \sim  \mathcal{P}(\alpha)$ be a random variable following a Poisson law of expectation $\alpha$. We consider the random graph $  \mathrm{G_{Poi}}(\alpha,p)$, which conditionally on $N$ is made of a classical Erd{\H{o}}s-R\'enyi $G(N,p)$ random graph (i.e. $N$ vertices where all ${N \choose 2}$ edges are independent and present with probability $p$) that we call the \textbf{core}, together with an infinite \textbf{stack} of vertices, which are all linked to every vertex of the $G(N,p)$ with probability $p$ independently. \emph{There is no edge  between vertices of the stack}. See Figure \ref{fig:ERP}. 

\paragraph{Markov property.} A step of \emph{exploration} in $  \mathrm{G}_ \mathrm{Poi}(\alpha,p)$ is the following: Fix a vertex $\rho$ of the stack and reveal its neighbors $ y_1, \dots , y_K$  with $K\geq 0$ inside the core. Then, see those vertices $y_1, \dots  , y_K$ as new vertices of the stack, in particular erase all possible connections between vertices of the stack. The key lemma is the following:
\begin{lemma}[Markov property of $ \mathrm{G_{Poi}}(\alpha,p)$] \label{lem:Markov} Let $K \geq 0$ be the number of neighbors in the core of a given vertex $\rho$ of the stack in $  \mathrm{G_{Poi}}(\alpha,p)$. Then $ K \sim \mathcal{P}( \alpha p)$ and conditionally on $K$, the graph made after removing $\rho$ and placing its  $K$ neighbors in the stack has law $ \mathrm{G_{Poi}}( \alpha (1-p),p)$.
\end{lemma}
\noindent \textbf{Proof.} Call $N' = N -K$ the remaining number of vertices after the revelation of the $K$ neighbors of $\rho$ in the $G(N,p)$ part. Then remark the following factorization:
 \begin{eqnarray*} \mathbb{P}(K=k \mbox{ and } N'=n) &=&  \mathrm{e}^{-\alpha} \frac{\alpha^{n+k}}{(n+k)!} \cdot {n+ k \choose k} p^{k} (1-p)^{n}\\  &=& \left( \mathrm{e}^{-\alpha(1-p)} \frac{(\alpha(1-p))^{n}}{n!}\right) \cdot \left( \mathrm{e}^{-\alpha p} \frac{(\alpha p)^{k}}{k!}\right).  \end{eqnarray*}
 The statement follows since conditionally on the status of the vertices (being in the stack, or in the remaining part), all possible edges are i.i.d.~present with probability $p$.
\qed 

\begin{figure}[!h]
       \begin{center}
       \includegraphics[width=12cm]{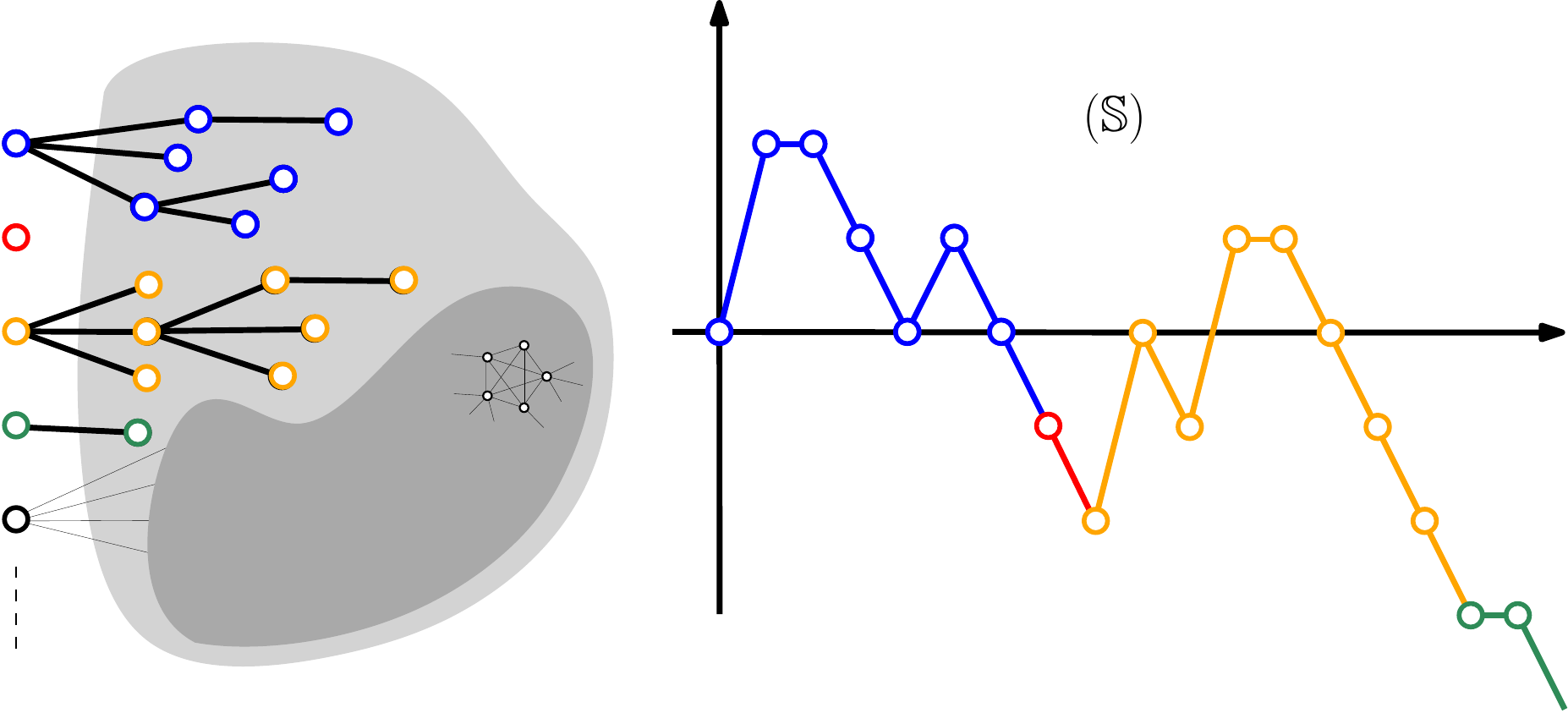}
       \caption{Lukasiewicz exploration of the graph $ \mathrm{G_{Poi}}(\alpha,p)$ obtained by concatenating the iterative number of neighbors $-1$. \label{fig:lukaER}}
       \end{center}
       \end{figure}

\paragraph{Lukasiewicz exploration.}
In particular, successive explorations in $ \mathrm{G_{Poi}}(\alpha,p)$ yields a sequence of \emph{independent} Poisson random variables with expectation $ \alpha p, \alpha p (1-p), \dots , \alpha p (1-p)^{k}, \dots$ whose total sum is just a Poisson variable of parameter $ \alpha p \sum_{i \geq 0} (1-p)^{i} = \alpha$, recovering the total number of vertices $N$ as expected. We shall now assume that iteratively, the vertices explored are placed on top of the stack and that we always explore the first vertex of the stack: we get the so-called Lukasiewicz exploration of the graph $ \mathrm{G_{Poi}}(\alpha,p)$, see Figure 	\ref{fig:lukaER}. We encode it in a process $ (\mathbb{S}_{k} : k \geq 0)$, the Lukasiewicz walk,  defined by $ \mathbb{S}_{0}=0$ and where $ \Delta \mathbb{S}_{k} = \mathbb{S}_{k}- \mathbb{S}_{k-1}$ is the number of neighbors discovered at step $k$ minus one. \emph{Each new minimal record of $ \mathbb{S}$ thus corresponds to the exploration of the connected component of a new vertex of the initial stack}. Using Lemma \ref{lem:Markov} we can write \emph{simultaneously} for all $k \geq 0$
  \begin{eqnarray}
   \mathbb{S}_{k} &= &(\mathcal{P}(\alpha p )-1) + ( \mathcal{P}(\alpha p (1-p))-1) + \dots + ( \mathcal{P}(\alpha p (1-p)^{k-1})-1)\nonumber \\ &=&   \mathcal{N}\left( \alpha p \cdot \sum_{i=0}^{k-1} (1-p)^{i} \right)-k = \mathcal{N}\left( \alpha (1-(1-p)^{k})\right) -k, \label{eq:lukapoisson} \end{eqnarray}
    where all the Poisson random variables written above are independent and where $ (\mathcal{N}(t): t \geq 0)$ is a standard unit-rate Poisson counting process on $ \mathbb{R}_{+}$. We shall only use the following  standard estimates on the Poisson counting  process 
      \begin{eqnarray} \label{eq:lawll}
      \frac{ \mathcal{N}(t)}{t} \xrightarrow[t\to\infty]{a.s.} 1, \quad \mbox{ and } \quad \frac{ ( \mathcal{N}(tn)-tn)}{ \sqrt{n}} \xrightarrow[n\to\infty]{(d)} (B_{t} : t \geq 0),  \end{eqnarray} where $(B_{t} : t \geq 0)$ is a standard linear Brownian motion. 
      
\section{Phase transition for the giant and Aldous' critical limit}
Let us use the Lukasiewicz exploration of the Poissonized version of the Erd{\H{o}}s--R\'enyi random graph to give a straightforward proof of the well-known phase transition for the size of the largest connected component.

\subsection{Existence of the giant component}
Fix $c>0$. Let $\alpha =n$ and $p \equiv p_n = \frac{c}{n}$ and denote by $ \mathbb{S}^{{(n)}}$ the resulting Lukasiewicz walk to emphasize the dependence in $n$. Since we have $(1- \frac{c}{n})^{[nt]} \to \mathrm{e}^{-ct}$  as $n \to \infty$ uniformly over compact time intervals, using \eqref{eq:lukapoisson} and the law of large numbers \eqref{eq:lawll} we immediately deduce:
\begin{proposition}[Fluid limit]\label{prop:fluid}We have the following convergence in probability for the uniform norm over every compact of $ \mathbb{R}_+$:
$$ \left(n^{-1} \cdot {\mathbb{S}^{(n)}_{[nt]}}\right)_{t \geq 0} \xrightarrow[n\to\infty]{( \mathbb{P})} \left( 1- \mathrm{e}^{{-ct}}-t\right)_{ t \geq 0}.$$
\end{proposition}
  \begin{figure}[!h]
  \begin{center}
  \includegraphics[width=10cm]{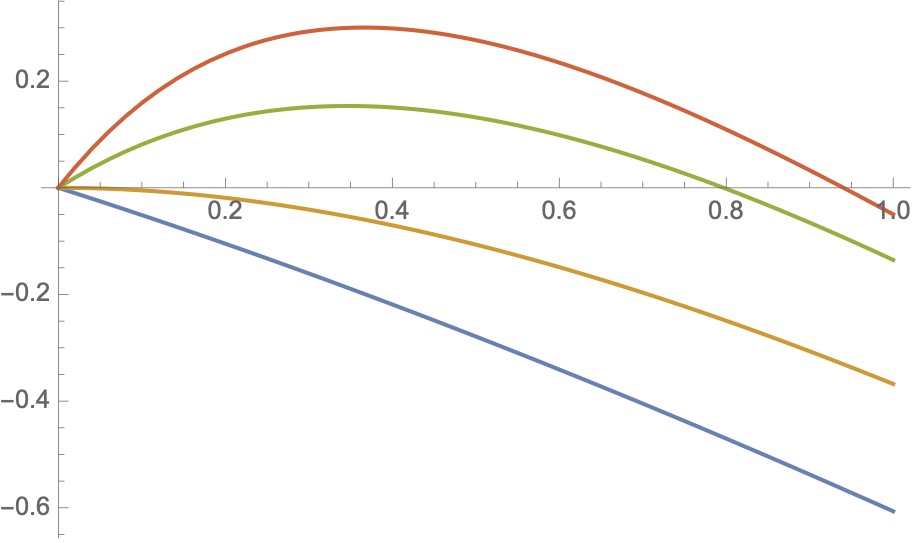}
  \caption{Graphs of the functions $( 1- \mathrm{e}^{{-ct}}-t)_{t \geq 0}$ for different of values of $c$: in blue $c=1/2$, in orange $c=1$, in green $c=2$ and in red $c=3$.}
  \end{center}
  \end{figure}
  
We write $X_{n} = o_{ \mathbb{P}}(n)$ if the random variable satisfies $X_{n}/n \to 0$ in probability.
 \begin{corollary}[Phase transition for $ \mathrm{G_{Poi}}(n, \frac{c}{n})$] \label{cor:giant} If $c <1$ then asymptotically all connected components of the core of $\mathrm{G_{Poi}}(n, \frac{c}{n})$ are $o_{ \mathbb{P}}(n)$, whereas if $ c >1$ it contains a unique giant component of size $ \beta(c)n + o_ \mathbb{P}(n)$ where $ \beta(c)$ is the first positive root of $ 1- \mathrm{e}^{{-ct}}-t =0$, all  others component sizes are $o_{ \mathbb{P}}(n)$.
 \end{corollary}
 \noindent \textbf{Proof.} The size of the connected components in $\mathrm{G_{Poi}}(n, \frac{c}{n})$ are given by the lengths of the excursions of $ \mathbb{S}^{(n)}$ above its running infimum process $ \underline{ \mathbb{S}}^{(n)}_k := \inf \{ \mathbb{S}^{(n)}_j : 0 \leq j \leq k\}$. When  $c>1$, fix $ \varepsilon \in (0, \beta(c)/2)$ and let us consider $$I^{(n)} = \inf\left\{ i \leq  \varepsilon n :  \mathbb{S}^{(n)}_i = \inf_{1 \leq j \leq \varepsilon n} \mathbb{S}^{(n)}_j\right\} \quad \mbox{ and } J^{(n)} = \inf\left\{ i \geq I^{(n)} : \mathbb{S}^{(n)}_i = \mathbb{S}^{(n)}_{I^{(n)}}-1\right\},$$ so that $ \mathbb{S}^{(n)}$ performs an excursion above is running infimum over the time interval $[I^{(n)}, J^{(n)}]$. By the property of the Lukasiewicz exploration, this excursion coincides with the exploration of a connected component $ \mathcal{C}^{{(n)}}$ of the stack of size $J^{(n)}-I^{(n)}$ and the above proposition entails that $I^{(n)}/n \to 0$ and $J^{(n)}/n \to \beta(c)$ in probability as $n \to \infty$. In particular we have $ | \mathcal{C}^{(n)}| = \beta(c) n + o_{ \mathbb{P}}(n)$ as desired. This component is further split into $\Delta \mathbb{S}_{I^{(n)}}+1$ components inside the core: Consider then
 $$ \tilde{I}^{(n)} = \inf\left\{ I^{(n)}+1 \leq i \leq  \varepsilon n :  \mathbb{S}^{(n)}_i = \inf_{I^{(n)}+1 \leq j \leq \varepsilon n} \mathbb{S}^{(n)}_j\right\} \quad \mbox{ and } \tilde{J}^{(n)} = \inf\left\{ i \geq \tilde{I}^{(n)} : \mathbb{S}^{(n)}_i = \mathbb{S}^{(n)}_{\tilde{I}^{(n)}}-1\right\}, $$ then the excursion time $[\tilde{I}^{(n)}, \tilde{J}^{(n)}]$ corresponds to the exploration of a subcomponent of $ \mathcal{C}^{(n)}$ inside the core (after removing the vertex of the stack it contains), and similarly we deduce from the above that $\tilde{I}^{(n)}/n \to 0$ and $\tilde{J}^{(n)}/n \to \beta(c)$ in probability as $n \to \infty$, thus showing the existence of the giant component inside the core.  The control of the other excursions above $ \mathbb{S}^{(n)}$ (and of all excursions in the case $c <1$) is done similarly using the convergence of Proposition \ref{prop:fluid}. A slight difficulty is that that the convergence over every compact of $ \mathbb{R}_+$ is not sufficient to prevent from unexpected deviations of the process $ \mathbb{S}^{(n)}$ at large values $k_{n}$ such that $k_n/n \to \infty$. But notice that the number of vertices left in the core after the first  $[xn]$ steps of exploration has law $ \mathcal{P}( n (1-c/n)^{[xn]})$ whose expectation is $\sim n \mathrm{e}^{-cx}$ as $n \to \infty$. By choosing $x$ large, we can ensure that with high probability all components of size comparable to $n$ are described by the $[xn]$ steps of $ \mathbb{S}^{(n)}$.  \qed 

  \paragraph{Back to the $G(n,p)$ model.} \label{paragraph:back}The analogous statements in the case of $ G(n,p)$ can be deduced from the above results. The key idea for the depoissonization being that $  \mathcal{P}(n)$ is concentrated around the value $n$ with $ \sqrt{n}$ fluctuations. Indeed, if we let  $N_{-} \sim \mathcal{P}(n -   n^{7/12})$ and $ N_{+} \sim \mathcal{P}( n +  n^{7/12})$ then we have the natural inclusions 
  \begin{eqnarray} \label{eq:sandwhich} G(N_{-},p) \subset G(n,p) \subset G(N_{+},p),  \end{eqnarray} which hold with high probability when $n \to \infty$ because $ \mathbb{P}(\mathcal{P}( n - n^{7/12}) \leq n) \to 1$ and $ \mathbb{P}(\mathcal{P}( n +n^{7/12}) \geq  n) \to 1$ as $n \to \infty$ (notice that we only coupled the cores of $ \mathrm{G_{Poi}}$ and not their stacks). By Corollary \ref{cor:giant}, with high probability, both random graphs on the left-hand side and right-hand side of the last display have a unique component of size $ \approx \beta(c)n$, all others being of size negligible in front of $n$. A moment of thought enables to deduce that the same is true for the graph $G(n,p)$  sandwhiched between those two, and this  is a classical result of Gilbert--Erd{\H{o}}s--R\'enyi \cite{erdds1959random,erdHos1960evolution}. 
 
\subsection{Refined estimates}
Let us turn to refined estimates on the clusters size still in the case $ \alpha=n$ and $p\equiv p_n =\frac{c}{n}$ for $ c >0$. 
\paragraph{CLT for the giant.} In this regime, for fixed $k\geq 0$ we have $n(1-(1- \frac{c}{n})^k) \to  c k$ as $n \to \infty$ so using \eqref{eq:lukapoisson} we deduce the following convergence in distribution 
$$ (\mathbb{S}^{(n)}_{k} : k \geq 0) \xrightarrow[n\to\infty]{} \left(\mathcal{N}( c \cdot k) -k : k \geq 0\right).$$ When $c >1$, using \eqref{eq:lawll} we easily deduce the convergence in law of the time $I^{(n)}$ defined in the course of the proof of Corollary \ref{cor:giant} towards $ \mathcal{I} := \mathrm{argmin}\{ \mathcal{N}(c \cdot k)-k : k \geq 0\}$. In particular, the number of connected components of $ \mathrm{G_{Poi}}$ explored before finding the giant converges to $ K_c = \mathcal{N}(c \cdot  \mathcal{I})- \mathcal{I}$ which by a simple application of Markov property is a geometric random variable (whose success parameter is a posteriori $\beta(c)$ by Corollary \ref{cor:giant}). Similarly, the time $J^{(n)}$ can be further estimated. Put $k = [\beta(c) n + x \sqrt{n}]$ for $x \in \mathbb{R}$ in a compact interval using \eqref{eq:lawll} and notice that as $n \to \infty$ we have 
$$ \mathbb{S}^{(n)}_{k} = \mathcal{N}\Big( \underbrace{n \left(1-(1- \frac{c}{n})^{[\beta(c)n + x \sqrt{n}]}\right)}_{ =  \beta(c) n  +  x c (1- \beta(c)) \sqrt{n}+ o(n^{1/2})}\Big)-k  \underset{ \eqref{eq:lawll}}{\approx}  \sqrt{n}\left(B_{ \beta(c)} + x(c(1- \beta(c))-1) + o_{ \mathbb{P}}(1)\right),$$ where $(B_t : t \geq 0)$ is the Brownian motion appearing in \eqref{eq:lawll}.  Combining those two arguments we easily deduce that   \begin{eqnarray} \label{eq:cltgiant} \frac{J^{(n)} - \beta(c)n}{ \sqrt{n}} \xrightarrow[n\to\infty]{(d)} \frac{B_{\beta{c}}}{1-c^*}, \quad \mbox{ where }c^{*} = c(1- \beta(c)),  \end{eqnarray} which together with the convergence of $I^{{(n)}}$ establishes the central limit theorem  for the size of the giant component in $ \mathrm{G_{Poi}}$. 
\begin{remark} In the case of the Erd{\H{o}}s--R\'enyi $G(n, \frac{c}{n})$ for $c>1$ the central limit theorem of the size of the giant makes appear a variance $ \beta(c)(1-\beta(c))/(1-c^{*})^{2}$, see  \cite{pittel1990tree,bollobas2012asymptotic}. The additional $(1- \beta(c))$ factor can easily be explained if we conditioned our model to have a core of size $ n + o( \sqrt{n})$. We suspect that it is possible to derive the fixed-size Erd{\H{o}}s--R\'enyi case from the above fact using soft arguments.\end{remark}

\paragraph{Critical case and Aldous's limit. } 
In the critical case $p = \frac{1}{n}$ we can give a analog of a result of Aldous in the case of the standard fixed-size Erd{\H{o}}s-R\'enyi \cite{aldous1997brownian} indicating that the size of the clusters in the near critical regime is of order $ n^{2/3}$:

\begin{proposition}[Near critical case] \label{prop:aldous} Fix $\lambda \in \mathbb{R}$. For $ p\equiv p_n = \frac{1}{n} + \frac{\lambda}{n^{{4/3}}}$ with $\lambda \in \mathbb{R}$, the Lukasiewicz walk of $ \mathrm{G_{Poi}}(n, p_n)$ satisfies 
$$  \left(n^{-1/3} \cdot {\mathbb{S}^{(n)}_{[n^{2/3}t]}}\right)_{t \geq 0} \xrightarrow[n\to\infty]{(d)} \left( B_{t} + \lambda t - \frac{t^{2}}{2}\right)_{ t \geq 0}.$$
\end{proposition}

\noindent \textbf{Proof.} Putting $k = [n^{2/3}t]$ for $t \in [0,A]$ in a compact time interval in the equation \eqref{eq:lukapoisson} yields to
  $$ \mathbb{S}^{(n)}_{[n^{2/3}t]} =   \mathcal{N}\Big( \underbrace{n \left(1-(1- \frac{1}{n}- \frac{\lambda}{n^{4/3}})^{[n^{2/3}t]}\right)}_{ =  tn^{2/3} + \lambda t n^{1/3} - \frac{t^{2}}{2} n^{1/3} + o(n^{1/3})}\Big) -[n^{2/3}t]   \underset{ \eqref{eq:lawll}}{\approx} n^{1/3} \left( B_{t} + \lambda t - \frac{t^{2}}{2}\right).\qquad \qed$$
  
  \begin{remark} We suspect it is possible to recover the result of Aldous on the fixed-size $G(n,p)$ using the above proposition and rather soft arguments. We however did not pursue this goal in this short note.
  \end{remark}

  \section{Connectedness}
  As another application of our Poissonization technique,  let us give a short proof of the sharp phase transition for connectedness in the fixed-size Erd{\H{o}}s--R\'enyi \cite{erdds1959random,erdHos1960evolution}:
  \begin{proposition} For $c \in \mathbb{R}$ we have 
  $$ \mathbb{P}\left( G\left(n, \frac{\log n +c}{n}\right) \mbox{ is connected}\right) \xrightarrow[n\to\infty]{} \mathrm{e}^{- \mathrm{e}^{-c}}.$$
  \end{proposition}
  \noindent \textbf{Proof.} Let $p\equiv p_n = \frac{\log n +c}{n}$. We shall first prove the convergence of the proposition when the number $N$ of vertices of the Erd{\H{o}}s--R\'enyi graph is random and distributed according to one plus a Poisson law of expectation $n$. Connectedness of this graph is equivalent to the fact that inside $ \mathrm{G_{Poi}}(n,p)$ all vertices of the stack, except the first one, have a trivial connected component. This is the case if and only if the Lukasiewicz walk $( \mathbb{S}^{(n)})$ starts with a (large) excursion and once it has reached level $-1$, it makes only jumps of $-1$. Using \eqref{eq:lukapoisson} and \eqref{eq:lawll}, it is easy to see that the first hitting time $\tau_{-1}^{(n)}$ of $-1$ by the process $ \mathbb{S}^{(n)}$ is concentrated around $n$ and more precisely using similar calculations as in \eqref{eq:cltgiant} shows that
 \begin{eqnarray} \label{eq:tau-1n} \frac{\tau_{-1}^{(n)}-n}{ \sqrt{n}} \xrightarrow[n\to\infty]{( d)} B_{1}.  \end{eqnarray}
Besides, since the increments of $ \mathbb{S}^{(n)}$ are Poisson and independent, by the Markov property of the exploration we have that
 \begin{eqnarray*} \mathbb{P}(  \rho \cup G(N,p) \mbox{ is connected}) &=&  
\mathbb{P}\left( \Delta \mathbb{S}^{(n)}_{i} =-1, \forall i \geq \tau_{-1}^{(n)} \right)=  \mathbb{E}\left[\mathbb{P}( \Delta \mathbb{S}^{(n)}_{i} =-1, \forall i \geq \tau_{-1}^{(n)} \mid \tau_{-1}^{(n)})\right] \\ 
&=& \mathbb{E}\left[ \mathbb{P}\left( \mathcal{P}(n (1-p_n)^{\tau_{-1}^{(n)}}) =0\right)\right] = \mathbb{E}\left[  \mathrm{exp}\left(- n\left (1- \frac{\log n +c }{n}\right)^{\tau_{-1}^{(n)}}\right)\right],  \end{eqnarray*} and by \eqref{eq:tau-1n} the random variable inside the expectation  converges in probability to $ \mathrm{e}^{-{ \mathrm{e}^{-c}}}$ since $\tau_{-1}^{{(n)}}  = n + O_{ \mathbb{P}}( \sqrt{n})$. The desired statement follows.
To come back to the fixed-size $G(n,p)$ model, notice that the function $ \phi(n,p) = \mathbb{P}( G(n,p) \mbox{ is connected})$  is increasing in $p$ for $n$ fixed, but the monotonicity in $n$ is not clear. However, the natural inclusion $ G(n,p) \subset G(n+1,p)$ enables to write 
 \begin{eqnarray*}  \phi(n+1,p)  &\geq&  \mathbb{P}(G(n,p) \mbox{ is connected and the $(n+1)th$ vertex is connected to one of the first $n$ vertices}) \\ &\geq & \phi(n,p) - \mathrm{P}( \mathrm{Bin}(n,p) =0) \geq \phi(n,p) -  \mathrm{e}^{-np}.   \end{eqnarray*} Recalling that $p_{n} = \frac{\log n +c}{n}$ we deduce that if $0\leq k_{n} = o(n)$  then we have $\phi(n+k_{n}, p_{n}) \geq \phi(n,p_{n}) + o(1)$. With the notation of \eqref{eq:sandwhich} this shows that
$$ \mathbb{P}( G(N_{-},p_{n}) \mbox{ is connected} )-o(1) \leq  \mathbb{P}(G(n,p_{n}) \mbox{ is connected}) \leq \mathbb{P}( G(N_{+},p_{n}) \mbox{ is connected} )+o(1),$$ and by sandwhiching,  the middle term does converge to $\mathrm{e}^{- \mathrm{e}^{-c}}$. \qed \medskip

\noindent  \textbf{Acknowledgments.} We thank Svante Janson for feedback on a first version of this note and to Yuval Peres and Bal\'azs R\'ath for pointing to us a few shortcomings which greatly helped improving the note.


\begin{thebibliography}{1}

\bibitem{aldous1997brownian}
{\sc D.~Aldous}, {\em {B}rownian excursions, critical random graphs and the
  multiplicative coalescent}, Ann. Probab.,  (1997), pp.~812--854.

\bibitem{bollobas2012asymptotic}
{\sc B.~Bollob{\'a}s and O.~Riordan}, {\em Asymptotic normality of the size of
  the giant component via a random walk}, Journal of Combinatorial Theory,
  Series B, 102 (2012), pp.~53--61.

\bibitem{erdds1959random}
{\sc P.~Erd{\H{o}}s and A.~R\'enyi}, {\em On random graphs i}, Publ. math.
  debrecen, 6 (1959), p.~18.

\bibitem{erdHos1960evolution}
{\sc P.~Erd{\H{o}}s and A.~R{\'e}nyi}, {\em On the evolution of random graphs},
  Publ. Math. Inst. Hung. Acad. Sci, 5 (1960), pp.~17--60.

\bibitem{pittel1990tree}
{\sc B.~Pittel}, {\em On tree census and the giant component in sparse random
  graphs}, Random Structures \& Algorithms, 1 (1990), pp.~311--342.

\end{thebibliography}
\end{document}